\documentclass[english]{smfart}

\usepackage[english]{babel}
\usepackage{smfthm}
\usepackage{bull}
\usepackage{mathrsfs}
\usepackage{amssymb,bm}
\usepackage{graphicx}
\usepackage{hyperref}
\usepackage{smfhyperref}
\usepackage[all]{xy}

\renewcommand{\bar}[1]{\overline{#1}}

\newcommand{\A}{\mathbb{A}}
\newcommand{\C}{\mathbb{C}}
\newcommand{\G}{\mathbb{G}}
\newcommand{\N}{\mathbb{N}}
\newcommand{\R}{\mathbb{R}}
\newcommand{\Z}{\mathbb{Z}}

\newcommand{\prim}{{\mathrm{P} \hspace{-0.3mm}}}
\newcommand{\lacc}{\{ \! \{}
\newcommand{\racc}{\} \! \}}

\NumberTheoremsAs{equation}

\author{Nicolas Martin}
\address{Centre de math\'ematiques Laurent Schwartz, \'Ecole polytechnique, Universit\'e Paris-Saclay, F-91128 Palaiseau cedex, France}
\email{nicolas.martin@polytechnique.edu}
\title{Middle multiplicative convolution and hypergeometric equations}

\begin{document}

\frontmatter

\begin{abstract}
Using a relation due to Katz linking up additive and multiplicative convolutions, we make explicit the behaviour of some Hodge invariants by middle
multiplicative convolution, following \cite{DS13} and \cite{Mar21} in the additive case. Moreover, the main theorem gives a new proof of a result of Fedorov computing the Hodge invariants of hypergeometric equations.
\end{abstract}


\maketitle

\mainmatter

\bigskip
\bigskip

The starting point of this article is a work of Dettweiler and Sabbah \cite{DS13} consisting in making explicit the behaviour of Hodge invariants by middle additive convolution by a Kummer module, motivated by the Katz algorithm \cite{Kat96}. In \cite{Mar21}, we developed this work without making the assumption of scalar monodromy at infinity assumed in the Katz algorithm and in \cite{DS13}, and we made precise the behaviour of nearby cycle local Hodge numerical data.

\bigskip

There exists a tricky link between middle additive convolution with a Kummer module and middle multiplicative convolution with a particular hypergeometric module, due to Katz \cite{Kat96} and detailed in Proposition \ref{relationkatz}. It allows us in {\S}\ref{2} to transpose the general results of \cite{Mar21} to the multiplicative context, after having recalled  in {\S}\ref{1} the necessary definitions.

\newpage

An application of these results is another way to prove a theorem due to Fedorov computing the Hodge invariants of hypergeometric equations \cite[Th.\,3]{Fed18} (detailed in \cite[Th.\,2.6]{Fou19}), very different but more direct, insofar as it uses the explicit behaviour of the Hodge invariants at infinity and 0.

\vspace{-0.1cm}

\section{Numerical Hodge data}

\label{1}
\numberwithin{equation}{section}

\vspace{0.1cm}

Let us begin by recalling the definition of local Hodge invariants introduced in \cite[{\S}2.2]{DS13}. Let $\Delta$ be a disc in the complex plane centered at 0 with coordinate $t$ and $(V,F^{\bullet}V,\nabla)$ be a polarizable variation of Hodge structure on $\Delta^*=\Delta \smallsetminus \{0\}$ over the field of complex numbers (for a complete review of this notion, see \cite[{\S}5.2]{Mar18}). We denote by $M$ the corresponding $\mathscr{D}_\Delta$-module minimal extension at $0$.

\vspace{0.3cm}

\noindent
\textbf{Nearby cycles.} For $a \in (-1,0]$ and $\lambda = \textnormal{exp}(-2i\pi a)$, the nearby cycle space at the origin $\psi_\lambda(M)$ is equipped with the nilpotent endomorphism $\mathrm{N} = t \partial_t-a$ and we have an induced Hodge filtration on $\psi_\lambda(M)$ such that $\mathrm{N} F^p\psi_\lambda(M) \subset F^{p-1}\psi_\lambda(M)$. The monodro\-my filtration induced by $\mathrm{N}$ enables us to define the spaces $\prim_\ell \psi_\lambda(M)$ of primitive vectors, equipped with a polarizable Hodge structure. The nearby cycle local Hodge numerical data are defined by

\vspace{-0.1cm}

$$\nu_{\lambda,\ell}^p(M):=h^p(\prim_\ell \psi_\lambda(M))=\dim \textnormal{gr}_F^p \prim_\ell \psi_\lambda(M),$$

\noindent
with the relation $\nu_{\lambda}^p(M):=h^p \psi_\lambda(M)=\sum\limits_{\ell \geq 0}\sum\limits_{k=0}^\ell \nu_{\lambda,\ell}^{p+k}(M)$. We set

\vspace{0.1cm}

$$\nu_{\lambda,\textnormal{prim}}^p(M):=\sum\limits_{\ell \geq 0} \nu_{\lambda,\ell}^p(M) \ \textnormal{ and } \ \nu_{\lambda,\textnormal{coprim}}^p(M):=\sum\limits_{\ell \geq 0} \nu_{\lambda,\ell}^{p+\ell}(M).$$

\vspace{0.1cm}

\noindent
\textbf{Vanishing cycles.} For $\lambda \neq 1$, the vanishing cycle space at the origin is given by $\phi_\lambda(M)=\psi_\lambda(M)$ and is equipped with $\textrm{N}$ and $F^p$ as before. For $\lambda=1$, the Hodge filtration on $\phi_1(M)$ is such that $F^p \prim_\ell \phi_1(M)=\mathrm{N}(F^p \prim_{\ell+1} \psi_1(M))$. Similarly to nearby cycles, the vanishing cycle local Hodge numerical data is defined by $\mu_{\lambda,\ell}^p(M):=h^p(\prim_\ell \phi_\lambda(M))=\dim \textnormal{gr}_F^p \prim_\ell \phi_\lambda(M)$.

\vspace{0.3cm}

Now let us leave the local point of view, and let $\bm{x}=\{x_1,...,x_r\}$ denote a set of points of $\G_m=\C^*$, $x_0=0$, $\mathscr{D}=\mathscr{D}_{\G_m}=\C[t,t^{-1}]\langle \partial_t \rangle$ and $i$ the inclusion $\G_m \!\! \smallsetminus \! \bm{x} \hookrightarrow \mathbb{P}^1$. Let $(V,F^\bullet V,\nabla)$ be a complex polarizable variation of Hodge structure on $\G_m \! \smallsetminus \bm{x}$ and $M$ be the $\mathscr{D}$-module minimal extension at points of~$\bm{x}$. Define $\mathscr{M}^{\textnormal{min}}$ to be the $\mathscr{D}_{\mathbb{P}^1}$-module minimal extension of $M$ at 0 and infinity.

\vspace{0.3cm}

\noindent
\textbf{Degrees $\delta^p$.} The Deligne extension $V^0$ of $(V,\nabla)$ on $\mathbb{P}^1$ is contained in $M$, and is endowed with the filtration $i_* F^p V \cap V^0$. We set
$$\delta^p(M)=\deg \textnormal{gr}_F^p V^0=\deg \, \frac{i_* F^p V \cap V^0}{i_* F^{p+1} V \cap V^0} \: .$$

\section{Middle multiplicative convolution with $H_{0,\gamma_0}$}

\label{2}
\numberwithin{equation}{section}

\vspace{0.3cm}

Let us fix $\gamma \in (0,1]$ and set $\lambda=\exp(-2i\pi \gamma)$. The Kummer module $\mathscr{L}_{\lambda}$ is defined by $\mathscr{L}_{\lambda}=\mathscr{D}/\mathscr{D} \cdot (t\partial_t - \gamma)$ where $\mathscr{D}=\mathscr{D}_{\G_m}$, and the middle additive convolution functor with $\mathscr{L}_{\lambda}$ is denoted by $\textnormal{MC}_{\lambda}$. Similarly to the middle additive convolution, the middle multiplicative convolution is defined by $M *_{\textnormal{mid} \times} N=\textnormal{Im}[\pi_\dagger(M \boxtimes N) \rightarrow \pi_+(M \boxtimes N)]$ where $M,N$ are holonomic $\mathscr{D}$-modules on $\G_m$, $\pi:\G_m \times \G_m \rightarrow \G_m$ is the product map, $\pi_+$ is the direct image functor, $\pi_\dagger:=\bm{D}\pi_+ \bm{D}$ is the adjoint by duality of $\pi_+$, and $M \boxtimes N$ denotes the external product of $M$ and $N$. See \cite[\S1.1]{DS13} for a quick review of middle convolution for holonomic modules on the affine line.

\medskip

Let us define $H_{0,\gamma}$ as the hypergeometric module $\mathscr{D}/\mathscr{D} \cdot (t\partial_t - t(t\partial_t-\gamma))$, whose restriction to $\G_m \smallsetminus \{1\}$ underlies a rank one local system with the following monodromies : 1 at 0, $\exp(-2i\pi\gamma)$ at 1, $\exp(2i\pi\gamma)$ at $\infty$. The next proposition links up additive and multiplicative convolutions and is due to Katz \cite[Lemma\,2.13.1]{Kat96}, and adapted here to the point of view of $\mathscr{D}$-modules (see \cite[Prop.\,2.8.1]{Mar18} for further details):

\vspace{0.1cm}

\begin{prop}\label{relationkatz}
Let us denote by $j:\G_{m} \hookrightarrow \A^1$ the inclusion. For every holonomic $\mathscr{D}$-module $M$, we have the following formula for $\gamma \neq 1$ :

\vspace{-0.15cm}

$$M *_{\textnormal{mid} \times} H_{0,\gamma} = j^+(\textnormal{MC}_\lambda(j_{\dagger +}(M \otimes \mathscr{L}_{\bar{\lambda}}))).$$
\end{prop}

\vspace{-0.1cm}

\begin{enonce}{Assumption}\label{assumption}\em
In everything that follows, we fix $\gamma_0 \in (0,1)$ and set $\lambda_0=\exp(-2i\pi \gamma_0)$. If we assume that $M$ is an irreducible regular holonomic $\mathscr{D}$-module, not isomorphic to $\mathscr{L}_{\lambda_0}$ and not supported on a point, then $j_{\dagger +}(M \otimes \mathscr{L}_{\bar{\lambda_0}})$ satisfies Assumption 1.2.2(1) of \cite{DS13} and we can apply to it the results of \cite{DS13} and \cite{Mar21}. Therefore, we make this assumption in what follows.\em
\end{enonce}

\vspace{0.1cm}

The following proposition gives the behaviour of vanishing cycle local Hodge numerical data by middle convolution with $H_{0,\gamma_0}$:

\vspace{0.1cm}

\begin{prop} \label{mup}
For all $i \in \{1,...,r\}$, we have:
\begin{equation*}
\mu_{x_i,\lambda,\ell}^p(M *_{\textnormal{mid} \times} H_{0,\gamma_0}) = \left\{
\begin{aligned}
  \mu_{x_i,\lambda/\lambda_0,\ell}^p(M) \quad &\textnormal{if } \gamma \in (0,\gamma_0] \\
	\mu_{x_i,\lambda/\lambda_0,\ell}^{p-1}(M) \quad &\textnormal{if } \gamma \in (\gamma_0,1].
\end{aligned}
\right.
\end{equation*}
\end{prop}

\vspace{0.15cm}

\begin{proof} For $i \in \{1,...,r\}$, Proposition \ref{relationkatz} gives
$$\mu_{x_i,\lambda,\ell}^p(M *_{\textnormal{mid} \times} H_{0,\gamma_0}) = \mu_{x_i,\lambda,\ell}^p(\textnormal{MC}_{\lambda_0}(j_{\dagger +}(M \otimes \mathscr{L}_{\bar{\lambda_0}}))).$$

\vspace{0.15cm}

\noindent
According to Assumption \ref{assumption}, we know that $j_{\dagger +}(M \otimes \mathscr{L}_{\bar{\lambda_0}})$ satisfies Assumption 1.2.2(1) of \cite{DS13}, then we can apply \cite[Th.\,3.1.2(2)]{DS13} (the part (2) does not need the hypothesis of scalar monodromy at infinity) and get

\vspace{-0.05cm}

\begin{equation*}
\mu_{x_i,\lambda,\ell}^p(M *_{\textnormal{mid} \times} H_{0,\gamma_0}) = \left\{
\begin{aligned}
  \mu_{x_i,\lambda/\lambda_0,\ell}^p(M \otimes \mathscr{L}_{\bar{\lambda_0}}) \quad &\textnormal{if } \gamma \in (0,\gamma_0] \\
	\mu_{x_i,\lambda/\lambda_0,\ell}^{p-1}(M \otimes \mathscr{L}_{\bar{\lambda_0}}) \quad &\textnormal{if } \gamma \in (\gamma_0,1].
\end{aligned}
\right.
\end{equation*}

\vspace{0.15cm}

\noindent
As $\mathscr{L}_{\bar{\lambda_0}}$ has trivial monodromy around $x_i \neq 0$, we have $\mu_{x_i,\lambda/\lambda_0,\ell}^p(M \otimes \mathscr{L}_{\bar{\lambda_0}})=\mu_{x_i,\lambda/\lambda_0,\ell}^p(M)$ and it is possible to conclude the proof.
\end{proof}

\vspace{0.15cm}

Concerning nearby cycle local Hodge numerical data at infinity, Theorem 1 of \cite{Mar21} complements \cite[Th.\,3.1.2]{DS13} without assuming that the monodromy at infinity is scalar. Combined with Proposition \ref{relationkatz} and \cite[2.2.13]{DS13}, we directly get the following proposition:

\vspace{0.15cm}

\begin{prop} \label{nuinfty}
We have the following data:
\vspace{0.2cm}
$$\nu_{\infty,\lambda,\ell}^p(M *_{\textnormal{mid} \times} H_{0,\gamma_0})=\left\{
\begin{array}{cl}
  \nu_{\infty,\lambda,\ell}^{p-1}(M) \quad &\textnormal{if } \gamma \in (0,1-\gamma_0) \\[1.5mm]
  \nu_{\infty,\lambda,\ell}^p(M) \quad &\textnormal{if } \gamma \in (1-\gamma_0,1) \\[1.5mm]
	\nu_{\infty,1,\ell+1}^{p}(M) \quad &\textnormal{if } \lambda=1 \\[1.5mm]
  \nu_{\infty,\bar{\lambda_0},\ell-1}^{p-1}(M) \quad &\textnormal{if } \lambda=\bar{\lambda_0}, \ \ell \geq 1.
\end{array}
\right.$$
\end{prop}

\vspace{-0.1cm}

\begin{rema}
We also have an explicit but more complicated formula for $\nu_{\infty,\bar{\lambda_0},0}^{p}(M *_{\textnormal{mid} \times} H_{0,\gamma_0})$, given and proved in \cite[Prop.\,6.4.3]{Mar18}.
\end{rema}

\vspace{0.15cm}

Similarly to Proposition \ref{nuinfty}, combining \cite[Th.\,3.1.2(2)]{DS13}, Proposition \ref{relationkatz} and \cite[2.2.14]{DS13}, the nearby cycle local Hodge numerical data at $0$ are given by the following proposition:

\vspace{0.15cm}

\begin{prop}[\protect{\cite[Prop.\,6.4.5]{Mar18}}]\label{nu0}
We have the following data:
\vspace{0.15cm}
$$\nu_{0,\lambda,\ell}^p(M *_{\textnormal{mid} \times} H_{0,\gamma_0}) = \left\{
\begin{array}{cl}
  \nu_{0,\lambda,\ell}^p(M) \quad &\textnormal{if } \gamma \in (0,\gamma_0) \\[1.5mm]
	\nu_{0,\lambda,\ell}^{p-1}(M) \quad &\textnormal{if } \gamma \in (\gamma_0,1) \\[1.5mm]
  \nu_{0,\lambda_0,\ell+1}^p(M) \quad &\textnormal{if } \lambda = \lambda_0 \\[1.5mm]
	\nu_{0,1,\ell-1}^{p-1}(M) \quad &\textnormal{if } \lambda = 1, \ \ell \geq 1 \\[1.5mm]
	h^p H^1(\mathbb{P}^1,\textnormal{DR}\mathscr{M}^{\textnormal{min}}) \quad &\textnormal{if } \lambda = 1, \ \ell = 0.
\end{array}
\right.$$
\end{prop}

\vspace{-0.1cm}

\begin{rema}
Summing the nearby cycle local Hodge numerical data, we deduce an explicit formula for Hodge numbers:

\vspace{-0.4cm}

\begin{multline*}
h^p(M *_{\textnormal{mid} \times} H_{0,\gamma_0}) = h^p(M) + \nu_{0,1,\textnormal{prim}}^{p-1}(M) - \nu_{0,\lambda_0,\textnormal{prim}}^{p-1}(M) \\
+ h^p H^1(\mathbb{P}^1,\textnormal{DR}\mathscr{M}^{\textnormal{min}}) + \sum_{\gamma \in [\gamma_0,1)} (\nu_{0,\lambda}^{p-1}(M)-\nu_{0,\lambda}^{p}(M)).
\end{multline*}

\end{rema}

\vspace{0.2cm}

To finish this study of the behaviour of Hodge invariants by middle multipli\-cative convolution with $H_{0,\gamma_0}$, let us make explicit the degrees $\delta^p$ defined in {\S}\ref{1}:\\

\vspace{-0.3cm}

\begin{prop}\label{deltap}
The degrees $\delta^p$ are given by:

\vspace{-0.4cm}

\begin{multline*}
\delta^p(M *_{\textnormal{mid} \times} H_{0,\gamma_0}) = \delta^p(M) + \sum_{\gamma \in [\gamma_0,1)} (\nu_{0,\lambda}^p(M)-\nu_{0,\lambda}^{p-1}(M)) + \nu_{0,\lambda_0,\textnormal{prim}}^{p-1}(M) \\
- \sum_{i=1}^r \left(\mu_{x_i,1}^p(M) + \sum_{\gamma \in (0,1-\gamma_0)} \mu_{x_i,\lambda}^{p-1}(M) \right).
\end{multline*}

\end{prop}

\vspace{-0.15cm}

\begin{proof} Applying \cite[Prop.\,2.3.2]{DS13} and \cite[2.2.13]{DS13}, we have

\begin{equation}\label{29}
\delta^p(M \otimes \mathscr{L}_{\bar{\lambda_0}})=\delta^p(M)- h^p(M)+ \!\! \sum_{\gamma \in [\gamma_0,1)} \!\! \nu_{0,\lambda}^p(M)+ \!\! \sum_{\gamma \in [1-\gamma_0,1)} \!\! \nu_{\infty,\lambda}^p(M)
\end{equation}

\begin{equation}\label{210}
\sum_{\gamma \in [\gamma_0,1)} \nu_{\infty,\lambda}^{p}(M \otimes \mathscr{L}_{\bar{\lambda_0}}) = \sum_{\gamma \in [\gamma_0,1)} \!\! \nu_{\infty,\lambda\bar{\lambda_0}}^{p}(M) = \sum_{\gamma \in [0,1-\gamma_0)} \nu_{\infty,\lambda}^{p}(M)
\end{equation}

\begin{equation}\label{211}
\sum_{\gamma \in (0,1-\gamma_0)} \mu_{0,\lambda}^{p-1}(M \otimes \mathscr{L}_{\bar{\lambda_0}}) = \sum_{\gamma \in (\gamma_0,1)} \nu_{0,\lambda}^{p-1}(M)
\end{equation}

\vspace{-0.1cm}

\begin{align}\label{212}
\mu_{0,1}^p(M \otimes \mathscr{L}_{\bar{\lambda_0}}) &= \nu_{0,1}^{p-1}(M \otimes \mathscr{L}_{\bar{\lambda_0}}) - \nu_{0,1,\textnormal{prim}}^{p-1}(M \otimes \mathscr{L}_{\bar{\lambda_0}}) \\
&= \nu_{0,\lambda_0}^{p-1}(M) - \nu_{0,\lambda_0,\textnormal{prim}}^{p-1}(M). \notag
\end{align}

\bigskip

\noindent
Summing \eqref{29} and \eqref{210}, we deduce from \cite[2.2.2($\ast\ast$)]{DS13} :

\vspace{-0.2cm}

\begin{equation}\label{213}
\delta^p(M \otimes \mathscr{L}_{\bar{\lambda_0}})+\sum_{\gamma \in [\gamma_0,1)} \nu_{\infty,\lambda}^{p}(M \otimes \mathscr{L}_{\bar{\lambda_0}})=\delta^p(M)+ \!\! \sum_{\gamma \in [\gamma_0,1)} \!\! \nu_{0,\lambda}^p(M).
\end{equation}

\vspace{0.2cm}

\noindent
According to Proposition \ref{relationkatz} and Theorem 3 of \cite{Mar21}, we have

\vspace{-0.3cm}

\begin{multline*}
\delta^p(M *_{\textnormal{mid} \times} H_{0,\gamma_0}) = \delta^p(M \otimes \mathscr{L}_{\bar{\lambda_0}}) + \sum_{\gamma \in [\gamma_0,1)} \nu_{\infty,\lambda}^{p}(M \otimes \mathscr{L}_{\bar{\lambda_0}}) \\
- \sum_{i=0}^r \left( \mu_{x_i,1}^p(M \otimes \mathscr{L}_{\bar{\lambda_0}}) + \sum_{\gamma \in (0,1-\gamma_0)} \mu_{x_i,\lambda}^{p-1}(M \otimes \mathscr{L}_{\bar{\lambda_0}}) \right).
\end{multline*}

\vspace{0.15cm}

\noindent
and, using \eqref{211}, \eqref{212} and \eqref{213}, we get the expected formula.
\end{proof}

\section{Fedorov's formula}

\label{3}

\vspace{0.2cm}

For any $\bm{\alpha},\bm{\beta} \in [0,1)^n$, the hypergeometric differential operator $\textnormal{Hyp}(\bm{\alpha},\bm{\beta})$ is defined by

\vspace{-0.4cm}

$$\textnormal{Hyp}(\bm{\alpha},\bm{\beta})=\prod_{i=1}^{n}(t \partial_t - \alpha_i)-t \prod_{j=1}^{n}(t \partial_t - \beta_j),$$

\noindent
and the corresponding hypergeometric module by $H_{\bm{\alpha},\bm{\beta}}:=\mathscr{D}/\mathscr{D} \cdot \textnormal{Hyp}(\bm{\alpha},\bm{\beta})$. These $\mathscr{D}$-modules are irreducibles if and only if $\alpha_i \neq \beta_j$ for all $i,j \in \{1,...,n\}$ \cite[Cor.\,3.2.1]{Kat90}. We assume in what follows that this condition is satisfied.

\vspace{0.3cm}

The leading term of the operator is $t^n(1-t)\partial_t^n$, then we have a connection on the trivial holomorphic bundle of rank $n$ on $\mathbb{P}^1 \smallsetminus \{0,1,\infty\}$. The three singularities are regular, and Theorem 3.5.4 of \cite{Kat90} shows that the corresponding local system on $\mathbb{P}^1 \smallsetminus \{0,1,\infty\}$ is physically rigid. In other words, the hypergeometric equation can be reconstructed, up to isomorphism, with the knowledge of its monodromies at $0$, $1$ and $\infty$, that was already remarked by Riemann in 1857. By \cite[Cor.\,8.1]{Sim90}, the restriction of $H_{\bm{\alpha},\bm{\beta}}$ to $\G_m \smallsetminus \{1\}$ underlies a complex polarizable variation of Hodge structure, unique up to a shift of the Hodge filtration \cite[Prop.\,1.13(i)]{Del87}. Let us make precise the monodromies of horizontal sections\footnote{\label{foot}Let us notice that Fedorov considers instead monodromies of solutions, see \cite[Prop.\,2.1]{Fed18}.} and the implications on the calculation of local Hodge invariants.

\vspace{0.2cm}

\noindent
\underline{\textbf{At $\infty$ :}} For $m \in \{1,...,n\}$, we set $\textnormal{mult}(\beta_m)= \# \{ j \in \{1,...,n\} \ | \ \beta_j=\beta_m\}$, $\ell_m(\bm{\beta})=\textnormal{mult}(\beta_m)-1$ and $\lambda_m=\exp(2i\pi\beta_m)$. The monodromy matrix at infinity is composed for each eigenvalue $\lambda_m$ of a unique Jordan block of size $\textnormal{mult}(\beta_m)$. We deduce that $\dim \prim_\ell \psi_{\infty,\lambda_m}(H_{\bm{\alpha},\bm{\beta}})=0$ except for $\ell = \ell_m(\bm{\beta})$ for which this quantity is equal to 1. The computation of $\nu_{\infty,\lambda_m,\ell}^p(H_{\bm{\alpha},\bm{\beta}})$ is reduced to finding the value of $p \in \Z$ for which this quantity for $\ell = \ell_m(\bm{\beta})$ is non zero (and equal to 1).

\vspace{0.2cm}

\noindent
\underline{\textbf{At $0$ :}} For $m \in \{1,...,n\}$, we set $\textnormal{mult}(\alpha_m)= \# \{ j \in \{1,...,n\} \ | \ \alpha_j=\alpha_m\}$, $\ell_m(\bm{\alpha})=\textnormal{mult}(\alpha_m)-1$ and $\mu_m=\exp(-2i\pi\alpha_m)$. The monodromy matrix at $0$ is composed for each eigenvalue $\mu_m$ of a unique Jordan block of size $\textnormal{mult}(\alpha_m)$. We deduce that $\dim \prim_\ell \psi_{0,\mu_m}(H_{\bm{\alpha},\bm{\beta}})=0$ except for $\ell = \ell_m(\bm{\alpha})$ for which this quantity is equal to 1. The computation of $\nu_{0,\mu_m,\ell}^p(H_{\bm{\alpha},\bm{\beta}})$ is reduced to finding the value of $p \in \Z$ for which this quantity for $\ell = \ell_m(\bm{\alpha})$ is non zero (and equal to 1).

\vspace{0.2cm}

\noindent
\underline{\textbf{At $1$ :}} Concerning the monodromy at $1$, we know from \cite[Prop.\,2.8]{BH89} that there are $n-1$ linearly independant eigenvectors associated to the eigenvalue 1 (see also \cite[Th.\,1.1]{Beu08}). That is, the monodromy at $1$ is a pseudoreflection. If we set $\gamma_s \in (0,1]$ such that $\gamma_s=\sum_{k=1}^n(\beta_k-\alpha_k)$ mod $\Z$, we deduce that $\lambda_s=\exp(-2i\pi\gamma_s)$ is also an eigenvalue of the monodromy at $1$, called the \textit{special eigenvalue}.

\vspace{0.1cm}

$\bullet$ If $\lambda_s \neq 1$, then the monodromy is diagonalizable. We have $\mu_{1,\lambda_s}(H_{\bm{\alpha},\bm{\beta}})=\nu_{1,\lambda_s}(H_{\bm{\alpha},\bm{\beta}})=1$, $\nu_{1,1}(H_{\bm{\alpha},\bm{\beta}})=n-1$ and $\mu_{1,1}(H_{\bm{\alpha},\bm{\beta}})=0$. The only thing left to be determined is the value of $p \in \Z$ for which $\mu_{1,\lambda_s,0}^p(H_{\bm{\alpha},\bm{\beta}})$ is non zero (and equal to 1).

\vspace{0.1cm}

$\bullet$ If $\lambda_s = 1$, then the monodromy is a transvection. We have $\nu_{1,1}(H_{\bm{\alpha},\bm{\beta}})=n$ and $\mu_{1,1}(H_{\bm{\alpha},\bm{\beta}})=1$. More precisely, $\mu_{1,1,\ell}(H_{\bm{\alpha},\bm{\beta}})=0$ except for $\ell=0$ for which this quantity is equal to 1. The only thing left to be determined is the value of $p \in \Z$ for which $\mu_{1,1,0}^p(H_{\bm{\alpha},\bm{\beta}})$ is non zero (and equal to 1).

\vspace{0.1cm}

\begin{defi}
Let us set $\alpha,\beta,\gamma \in [0,1)$. We say that the pair $(\alpha,\beta)$ \em is separated by \em $\gamma$ if $\exp(2i\pi\gamma)$ is in the open interval $(\exp(2i\pi\alpha),\exp(2i\pi\beta))$ of the oriented circle, a property that we denote by $\alpha \rightarrow \gamma \rightarrow \beta$. It means that either $0 \leq \alpha < \gamma < \beta < 1$, or $0 \leq \gamma < \beta < \alpha < 1$, or $0 \leq \beta < \alpha < \gamma < 1$.
\end{defi}

\vspace{-0.1cm}

\begin{rema}
It is the same notation as in the beginning of Chapter 4 of \cite{Fed18}, with the difference that $\alpha$, $\beta$ and $\gamma$ are not necessarily distinct (if they are not distinct, our property $\alpha \rightarrow \gamma \rightarrow \beta$ is not satisfied).
\end{rema}

\vspace{-0.1cm}

\begin{defi}
For $\bm{\alpha},\bm{\beta} \in [0,1)^n$ and $\gamma \in [0,1)$, we set
$$p(\bm{\alpha},\bm{\beta},\gamma) := \#\{k \ | \ \neg(\alpha_{k} \rightarrow \gamma \rightarrow \beta_{k}) \}=n-\# \, \{k \ | \ \alpha_{k} \rightarrow \gamma \rightarrow \beta_{k} \}.$$
Note that this quantity does not depend on the numbering of the $n$-tuple of pairs $((\alpha_1,\beta_1),...,(\alpha_n,\beta_n))$.
\end{defi}

\vspace{-0.1cm}

\begin{rema} According to \cite[Th.\,5.3.1]{Kat90}, given any partitions $\bm{\alpha} = \bm{\alpha'} \sqcup \bm{\alpha''}$ and $\bm{\beta} = \bm{\beta'} \sqcup \bm{\beta''}$ with $\#\bm{\alpha'}=\#\bm{\beta'}$, there exists a decomposition
$$H_{\bm{\alpha},\bm{\beta}} = H_{\bm{\alpha'},\bm{\beta'}} * H_{\bm{\alpha''},\bm{\beta''}} = H_{\bm{\alpha''},\bm{\beta''}} * H_{\bm{\alpha'},\bm{\beta'}}$$
where $*$ is the $*$-multiplicative convolution. In particular, for such partitions, the right-hand side is an object of $D_{\textnormal{hol}}^b(\mathscr{D})$ with cohomology in degree zero only. Furthermore, one can check that this formula also holds for the $!$-multiplicative convolution, hence for the middle multiplicative convolution. In the following, $*$ denotes any of these convolutions.

\vspace{0.1cm}

Moreover, for three-term partitions with equal corresponding size, the convolution $H_{\bm{\alpha},\bm{\beta}} * H_{\bm{\alpha'},\bm{\beta'}} * H_{\bm{\alpha''},\bm{\beta''}}$ is associative, as indicated in \cite[(5.1.7)]{Kat90}. As a consequence, given any permutation $\sigma \in \mathfrak{S}_n$, there exists a decomposition $H_{\bm{\alpha},\bm{\beta}} = H_{\alpha_1,\beta_{\sigma(1)}} * \cdots * H_{\alpha_n,\beta_{\sigma(n)}}$ which can be performed in a commutative and associative way.
\end{rema}

\vspace{-0.1cm}

\begin{rema}
1) Given a \textbf{fixed} decomposition into convolutions of hypergeometric modules of rank one, there exists a unique associated Hodge filtration by computing the convolutions one after the other, if we started from the trivial Hodge filtration for rank one : $F^p H_{\alpha,\beta}=H_{\alpha,\beta}$ for $p \leq 0$ and $F^p H_{\alpha,\beta}=0$ for $p \geq 1$.\\
2) By uniqueness of the Hodge filtration on $H_{\bm{\alpha},\bm{\beta}}$ up to a shift, we deduce that taking another decomposition will induce a shift in the filtration.
\end{rema}

\vspace{0.15cm}

For $a \in \R$, we denote by $\lacc a \racc$ the representative of $a \textnormal{ mod } \Z$ in $(0,1]$. We have $\lacc a \racc=a-\lceil a \rceil +1$.

\vspace{0.2cm}

\begin{theo}\label{fedorov}
Let us consider the decomposition $H_{\alpha_1,\beta_1} * \cdots * H_{\alpha_n,\beta_n}$ of $H_{\bm{\alpha},\bm{\beta}}$ into convolutions of hypergeometric modules of rank $1$. The variation of Hodge structure on $H_{\bm{\alpha},\bm{\beta}}$ induced by this decomposition satisfies:

\vspace{-0.3cm}

\begin{equation}\tag{a}
\nu_{0,\mu_m,\ell}^p(H_{\bm{\alpha},\bm{\beta}})=\left\{
\begin{array}{cl}
  \ 1 \quad &\textnormal{if } p=p(\bm{\alpha},\bm{\beta},\alpha_m) \textnormal{ and } \ell=\ell_m(\bm{\alpha}) \\[1mm]
  \ 0 \quad &\textnormal{otherwise}
\end{array}
\right.
\end{equation}

\vspace{-0.2cm}

\begin{equation}\tag{b}
\nu_{\infty,\lambda_m,\ell}^p(H_{\bm{\alpha},\bm{\beta}})=\left\{
\begin{array}{cl}
  \ 1 \quad &\textnormal{if } p=p(\bm{\alpha},\bm{\beta},\beta_m) \textnormal{ and } \ell=\ell_m(\bm{\beta}) \\[1mm]
  \ 0 \quad &\textnormal{otherwise}
\end{array}
\right.
\end{equation}

\vspace{-0.15cm}

\begin{equation}\tag{c}
\mu_{1,\lambda_s,\ell}^p(H_{\bm{\alpha},\bm{\beta}})=\left\{
\begin{array}{cl}
  \ 1 \ \ &\textnormal{if } p=n- \displaystyle{\left\lceil \sum_{k=1}^n \lacc \beta_k-\alpha_k \racc \right\rceil} \textnormal{ and } \ell=0 \\[1mm]
  \ 0 \ \ &\textnormal{otherwise.}
\end{array}
\right.
\end{equation}
\end{theo}

\vspace{0.1cm}

\begin{proof} By induction on $n \in \N^{*}$, that is, on the length of $\bm{\alpha}$ and $\bm{\beta}$. The theorem is satisfied for $n=1$. Let us set $n \geq 1$, $(\bm{\alpha},\bm{\beta})=((\alpha_0,...,\alpha_n),(\beta_0,...,\beta_n))$ two $(n+1)$-tuples such that $\alpha_i \neq \beta_j$ for all $i,j \in \{0,...,n\}$, and $m \in \{0,...,n\}$. We denote by $\{ \cdot \}$ the fractional part.

\vspace{0.3cm}

\noindent
\textbf{Formula (b).} Let us suppose that (b) is satisfied for all tuples of length $n$. Let us prove the formula for $\bm{\alpha}$ and $\bm{\beta}$ of length $n+1$.\\
(Case 1) Let us suppose that $\beta_m \neq \beta_0$. By physical rigidity and according to \cite[2.2.13]{DS13}, we have

\vspace{-0.1cm}

$$\nu_{\infty,\lambda_m,\ell}^p(H_{\bm{\alpha},\bm{\beta}})=\nu_{\infty,\lambda_m \! \exp(-2i\pi\alpha_0),\ell}^p(H_{\{\bm{\alpha}-\alpha_0\},\{\bm{\beta}-\alpha_0\}}),$$

\vspace{0.2cm}

\noindent
where $\{\bm{\alpha}-\alpha_0\}=(\{\alpha_0-\alpha_0\},\{\alpha_1-\alpha_0\},\hdots,\{\alpha_n-\alpha_0\})$. Then we have

\vspace{-0.1cm}

$$H_{\{\bm{\alpha}-\alpha_0\},\{\bm{\beta}-\alpha_0\}} = H_{\{\widehat{\bm{\alpha}_{0}}-\alpha_0\},\{\widehat{\bm{\beta}_0}-\alpha_0\}} * H_{0,\{\beta_0-\alpha_0\}},$$

\vspace{0.15cm}

\noindent
where $\widehat{\bm{\alpha}_{0}}$ is the tuple $\bm{\alpha}$ where we have removed $\alpha_0$, and $\widehat{\bm{\beta}_{0}}$ is the tuple $\bm{\beta}$ where we have removed $\beta_0$.

\vspace{0.1cm}

\noindent
Applying Proposition \ref{nuinfty}, we get

\vspace{-0.4cm}

$$\nu_{\infty,\lambda_m \! \exp(-2i\pi\alpha_0),\ell}^p(H_{\{\bm{\alpha}-\alpha_0\},\{\bm{\beta}-\alpha_0\}}) \! = \! \left\{
\begin{array}{cl}
  \nu_{\infty,\lambda_m \! \exp(-2i\pi\alpha_0),\ell}^{p-1} \! \left( H_{\{\widehat{\bm{\alpha}_{0}}-\alpha_0\},\{\widehat{\bm{\beta}_0}-\alpha_0\}} \right) \\[2mm]
	\textnormal{if } \{\beta_m-\alpha_0\} > \{\beta_0-\alpha_0\} \\[1mm]
  \nu_{\infty,\lambda_m \! \exp(-2i\pi\alpha_0),\ell}^{p} \! \left( H_{\{\widehat{\bm{\alpha}_{0}}-\alpha_0\},\{\widehat{\bm{\beta}_0}-\alpha_0\}} \right) \\[2mm]
	\textnormal{if } \{\beta_m-\alpha_0\} < \{\beta_0-\alpha_0\}.
\end{array}
\right.$$

\vspace{0.05cm}

\noindent
Applying \cite[2.2.13]{DS13} once again, we have

\vspace{-0.2cm}

$$\nu_{\infty,\lambda_m,\ell}^p(H_{\bm{\alpha},\bm{\beta}}) = \left\{
\begin{array}{cl}
  \nu_{\infty,\lambda_m,\ell}^{p} \left( H_{\widehat{\bm{\alpha}_0},\widehat{\bm{\beta}_0}} \right) \quad &\textnormal{if } \alpha_0 \rightarrow \beta_m \rightarrow \beta_0 \\
  \nu_{\infty,\lambda_m,\ell}^{p-1} \left( H_{\widehat{\bm{\alpha}_0},\widehat{\bm{\beta}_0}} \right) \quad &\textnormal{otherwise.}
\end{array}
\right.$$

\vspace{0.2cm}

\noindent
By the induction hypothesis, the left quantity is non zero if and only if $p=p(\bm{\alpha},\bm{\beta},\beta_m)$ and $\ell=\ell_m(\bm{\beta})=\ell_m(\widehat{\bm{\beta}_0})$.

\vspace{0.2cm}

\noindent
(Case 2) Let us suppose that $\beta_m = \beta_0$ and $\ell_0(\bm{\beta}) \geq 1$. Applying the same reasoning as before and using Proposition \ref{nuinfty} (case $\lambda=\bar{\lambda_0}$, $\ell \geq 1$), we get

\vspace{-0.2cm}

$$\nu_{\infty,\lambda_0,\ell}^p(H_{\bm{\alpha},\bm{\beta}})=\nu_{\infty,\lambda_0,\ell-1}^{p-1} \left( H_{\widehat{\bm{\alpha}_0},\widehat{\bm{\beta}_0}} \right) ,$$

\vspace{0.2cm}

\noindent
non zero if and only if $\ell=\ell_0(\bm{\beta})=\ell_0(\widehat{\bm{\beta}_0})+1$. In this case, we have $p(\bm{\alpha},\bm{\beta},\beta_0)=p(\widehat{\bm{\alpha}_0},\widehat{\bm{\beta}_0},\beta_0)+1$ because we do not have $\alpha_0 \rightarrow \beta_0 \rightarrow \beta_0$.

\vspace{0.2cm}

\noindent
(Case 3) Let us suppose that $\beta_m = \beta_0$ and $\ell_0(\bm{\beta}) = 0$, so we have $\beta_1 \neq \beta_0$. Applying the same reasoning as in Case 1, we get

\vspace{-0.1cm}

$$\nu_{\infty,\lambda_0,\ell}^p(H_{\bm{\alpha},\bm{\beta}}) = \left\{
\begin{array}{cl}
  \nu_{\infty,\lambda_0,\ell}^{p} \left( H_{\widehat{\bm{\alpha}_1},\widehat{\bm{\beta}_1}} \right) \quad &\textnormal{if } \{\beta_0-\alpha_1\} < \{\beta_1-\alpha_1\} \\
  \nu_{\infty,\lambda_0,\ell}^{p-1} \left( H_{\widehat{\bm{\alpha}_1},\widehat{\bm{\beta}_1}} \right) \quad &\textnormal{if } \{\beta_0-\alpha_1\} > \{\beta_1-\alpha_1\}
\end{array}
\right.$$
$$ \quad \quad \ \ = \left\{
\begin{array}{cl}
  \nu_{\infty,\lambda_0,\ell}^{p} \left( H_{\widehat{\bm{\alpha}_1},\widehat{\bm{\beta}_1}} \right) \quad &\textnormal{if } \alpha_1 \rightarrow \beta_0 \rightarrow \beta_1 \\
  \nu_{\infty,\lambda_0,\ell}^{p-1} \left( H_{\widehat{\bm{\alpha}_1},\widehat{\bm{\beta}_1}} \right) \quad &\textnormal{otherwise.}
\end{array}
\right.$$

\vspace{0.2cm}

\noindent
By the induction hypothesis, and as the order in which the convolutions are done does not matter, the left quantity is non zero if and only if $p=p(\bm{\alpha},\bm{\beta},\beta_0)$ and $\ell=\ell_0(\bm{\beta})=\ell_0(\widehat{\bm{\beta}_1})=0$.

\vspace{0.2cm}

\noindent
To conclude, Formula (b) is satisfied for the couple $(\bm{\alpha},\bm{\beta})$.

\vspace{0.3cm}

\noindent
\textbf{Formula (a).} Let us suppose that (a) is satisfied for all tuples of length $n$. Let us prove the formula for $\bm{\alpha}$ and $\bm{\beta}$ of length $n+1$.\\
(Case 1) Let us suppose that $\alpha_m \neq \alpha_0$. According to Proposition \ref{nu0} and \cite[2.2.13]{DS13}, and applying the same reasoning as in Case 1 of the proof of Formula (b), we have

$$\nu_{0,\mu_m,\ell}^p(H_{\bm{\alpha},\bm{\beta}}) = \left\{
\begin{array}{cl}
  \nu_{0,\mu_m,\ell}^{p} \left( H_{\widehat{\bm{\alpha}_0},\widehat{\bm{\beta}_0}} \right) \quad &\textnormal{if } \alpha_0 \rightarrow \alpha_m \rightarrow \beta_0 \\
  \nu_{0,\mu_m,\ell}^{p-1} \left( H_{\widehat{\bm{\alpha}_0},\widehat{\bm{\beta}_0}} \right) \quad &\textnormal{otherwise.}
\end{array}
\right.$$

\vspace{0.2cm}

\noindent
By the induction hypothesis, the left quantity is non zero if and only if $p=p(\bm{\alpha},\bm{\beta},\alpha_m)$ and $\ell=\ell_m(\bm{\alpha})=\ell_m(\widehat{\bm{\alpha}_0})$.

\vspace{0.2cm}

\noindent
(Case 2) Let us suppose that $\alpha_m = \alpha_0$ and $\ell_0(\bm{\alpha}) \geq 1$. Applying the same reasoning as before and using Proposition \ref{nu0} (case $\lambda=1$, $\ell \geq 1$), we get

\vspace{-0.1cm}

$$\nu_{0,\mu_0,\ell}^p(H_{\bm{\alpha},\bm{\beta}})=\nu_{0,\mu_0,\ell-1}^{p-1} \left( H_{\widehat{\bm{\alpha}_0},\widehat{\bm{\beta}_0}} \right) ,$$

\vspace{0.1cm}

\noindent
non zero if and only if $\ell=\ell_0(\bm{\alpha})=\ell_0(\widehat{\bm{\alpha}_0})+1$. In this case, we have $p(\bm{\alpha},\bm{\beta},\alpha_0) = p(\widehat{\bm{\alpha}_0},\widehat{\bm{\beta}_0},\alpha_0)+1$ because we do not have $\alpha_0 \rightarrow \alpha_0 \rightarrow \beta_0$.

\vspace{0.2cm}

\noindent
(Case 3) Let us suppose that $\alpha_m = \alpha_0$ and $\ell_0(\bm{\alpha}) = 0$, so we have $\alpha_1 \neq \alpha_0$. Applying the same reasoning as in Case 1, we get

\vspace{-0.2cm}

$$\nu_{0,\mu_0,\ell}^p(H_{\bm{\alpha},\bm{\beta}}) = \left\{
\begin{array}{cl}
  \nu_{0,\mu_0,\ell}^{p} \left( H_{\widehat{\bm{\alpha}_1},\widehat{\bm{\beta}_1}} \right) \quad &\textnormal{if } \alpha_1 \rightarrow \alpha_0 \rightarrow \beta_1 \\
  \nu_{0,\mu_0,\ell}^{p-1} \left( H_{\widehat{\bm{\alpha}_1},\widehat{\bm{\beta}_1}} \right) \quad &\textnormal{otherwise.}
\end{array}
\right.$$

\vspace{0.2cm}

\noindent
By the induction hypothesis, and as the order in which the convolutions are done does not matter, the left quantity is non zero if and only if $p=p(\bm{\alpha},\bm{\beta},\alpha_0)$ and $\ell=\ell_0(\bm{\alpha})=\ell_0(\widehat{\bm{\alpha}_1})=0$.

\vspace{0.2cm}

\noindent
To conclude, Formula (a) is satisfied for the couple $(\bm{\alpha},\bm{\beta})$.

\vspace{0.3cm}

\noindent
\textbf{Formula (c).} Let us suppose that Formula (c) is satisfied for all tuples of length $n$. Let us prove the formula for $\bm{\alpha}$ and $\bm{\beta}$ of length $n+1$. We denote by $\lambda_s$ the special eigenvalue of $H_{\bm{\alpha},\bm{\beta}}$, and by $\lambda'_s$ the special eigenvalue of $H_{\widehat{\bm{\alpha}_0},\widehat{\bm{\beta}_0}}$ and for $i \in \{1,...,n\}$ :
\begin{center}
$\gamma_i= \lacc \beta_i-\alpha_i \racc, \qquad \gamma_{\geq i}=\lacc \sum_{k \geq i}(\beta_k-\alpha_k) \racc = \lacc \sum_{k \geq i}\gamma_k \racc$
\end{center}
With these notations, we have $\gamma_{\geq 0}=\gamma_s$ and $\gamma_{\geq 1}=\gamma'_s$. Let us remark that for $\gamma_0,...,\gamma_n \in (0,1]$, we have the following relation :

\vspace{-0.25cm}

\begin{equation}
\label{eq}
\lceil \gamma_0 + \cdots + \gamma_n \rceil=\left\{
\begin{array}{cl}
  \lceil \gamma_1 + \cdots + \gamma_n \rceil \quad &\textnormal{if } \gamma_0+\gamma_{\geq 1} \leq 1 \\[1.5mm]
  \lceil \gamma_1 + \cdots + \gamma_n \rceil+1 \quad &\textnormal{if } \gamma_0+\gamma_{\geq 1} > 1.
\end{array}
\right.
\end{equation}

\begin{proof} Firstly $\gamma_0+\gamma_{\geq 1}=\gamma_0+\lacc\sum_{k \geq 1}\gamma_k\racc=\sum_{k \geq 0}\gamma_k-\lceil\sum_{k \geq 1}\gamma_k\rceil+1$, then\\

\vspace{-0.4cm}

\begin{center}
$\gamma_0+\gamma_{\geq 1} \leq 1 \ \Longleftrightarrow \ \sum_{k \geq 0}\gamma_k \leq \lceil\sum_{k \geq 1}\gamma_k\rceil \ \Longleftrightarrow \ \lceil\sum_{k \geq 0}\gamma_k\rceil=\lceil\sum_{k \geq 1}\gamma_k\rceil$
\end{center}
since $\gamma_0>0$. Similarly, $\gamma_0+\gamma_{\geq 1} > 1 \Longleftrightarrow \lceil\sum_{k \geq 0}\gamma_k\rceil=\lceil\sum_{k \geq 1}\gamma_k\rceil+1$.
\end{proof}

\vspace{0.2cm}

\noindent
Now, according to Proposition \ref{mup} and \cite[2.2.14]{DS13}, and applying the same reasoning as in the proof of Case 1 of Formula (b), we have

\vspace{-0.15cm}

$$\mu_{1,\lambda_s,\ell}^p(H_{\bm{\alpha},\bm{\beta}})=\left\{
\begin{array}{cl}
  \mu_{1,\lambda'_s,\ell}^{p} \left( H_{\widehat{\bm{\alpha}_0},\widehat{\bm{\beta}_0}} \right) \quad &\textnormal{if } \gamma_s \in (0,\gamma_0] \\[1.5mm]
  \mu_{1,\lambda'_s,\ell}^{p-1} \left( H_{\widehat{\bm{\alpha}_0},\widehat{\bm{\beta}_0}} \right) \quad &\textnormal{if } \gamma_s \in (\gamma_0,1].
\end{array}
\right.$$

\vspace{0.2cm}

\noindent
As $\gamma_s=\lacc \gamma_0+\gamma_{\geq1}\racc$, the condition $\gamma_s \in (0,\gamma_0]$ is equivalent to $\gamma_0+\gamma_{\geq1} > 1$. Likewise, the condition $\gamma_s \in (\gamma_0,1]$ is equivalent to $\gamma_0+\gamma_{\geq1} \leq 1$. By induction hypothesis, we deduce that

\vspace{-0.3cm}

\begin{align*}
\mu_{1,\lambda_s,\ell}^p(H_{\bm{\alpha},\bm{\beta}}) &= \left\{
\begin{array}{cl}
  \ 1 \ \ &\textnormal{if } p=n-\lceil\sum_{k \geq 1}\gamma_k\rceil \textnormal{ and } \gamma_0+\gamma_{\geq1} > 1  \textnormal{ and } \ell=0\\
	\ 1 \ \ &\textnormal{if } p=n-\lceil\sum_{k \geq 1}\gamma_k\rceil+1 \textnormal{ and } \gamma_0+\gamma_{\geq1} \leq 1 \textnormal{ and } \ell=0\\
  \ 0 \ \ &\textnormal{otherwise.}
\end{array}
\right.\\
&= \left\{
\begin{array}{cl}
  \ 1 \ \ &\textnormal{if } p=n+1-\lceil\sum_{k \geq 0}\gamma_k\rceil \textnormal{ and } \ell=0\\
  \ 0 \ \ &\textnormal{otherwise.}
\end{array}
\right.
\end{align*}

\vspace{0.1cm}

\noindent
according to \eqref{eq}. Then Formula (c) is satisfied for the couple $(\bm{\alpha},\bm{\beta})$.
\end{proof}

\vspace{0.2cm}

\noindent
\textbf{Link between Theorem \ref{fedorov} and Fedorov's formulas.} Formulas (a) and (b) of the previous theorem corresponds to Formulas (a) and (b) of Theorem 3 in \cite{Fed18}. However, this is not fully obvious in the sense that Fedorov considers in his article the space of solutions of the connection associated with the hypergeometric equation, while we consider the space of horizontal sections of the connection (see Footnote \ref{foot}). Let us begin by transposing Fedorov's formulas in terms of horizontal sections with the following lemma. Note that we do not necessarily assume that the tuples are ordered.

\vspace{0.1cm}

\begin{lemm} Parts $\textnormal{(a)}$ and $\textnormal{(b)}$ of \cite[Th.\,3]{Fed18} are equivalent to the follo\-wing statement:\\
The hypergeometric module $H_{\bm{\alpha},\bm{\beta}}$ is equipped with a polarizable variation of Hodge structures verifying, up to a shift, the following identities:

\vspace{-0.2cm}

\begin{equation}\tag{a}
\nu_{0,\mu_m,\ell}^p(H_{\bm{\alpha},\bm{\beta}})=\left\{
\begin{array}{cl}
  \ 1 \quad &\textnormal{if } p=\#\{j \ | \ \beta_j < \alpha_m \} - \#\{i \ | \ \alpha_i < \alpha_m \} \\[1.5mm]
	&\textnormal{and } \ell=\ell_m(\bm{\alpha}) \\[1.5mm]
  \ 0 \quad &\textnormal{otherwise.}
\end{array}
\right.
\end{equation}

\vspace{-0.05cm}

\begin{equation}\tag{b}
\nu_{\infty,\lambda_m,\ell}^p(H_{\bm{\alpha},\bm{\beta}})=\left\{
\begin{array}{cl}
  \ 1 \quad &\textnormal{if } p=\#\{j \ | \ \beta_j \leq \beta_m \} - \#\{i \ | \ \alpha_i < \beta_m \} \\[1.5mm]
	&\textnormal{and } \ell=\ell_m(\bm{\beta}) \\[1.5mm]
  \ 0 \quad &\textnormal{otherwise.}
\end{array}
\right.
\end{equation}

\end{lemm}

\vspace{0.1cm}

\begin{proof} The space of solutions and the space of horizontal sections are dual (see for example \cite[Cor.\,7.1.1]{Pha79}). If we denote by $*$ the dual, we have the relation $(\prim_\ell H)^* \simeq \textrm{N}^\ell \prim_\ell(H^*)$ as Hodge structures and then

\vspace{-0.05cm}

$$(\textnormal{gr}_F^p \prim_\ell H)^* \simeq \textnormal{gr}_F^{-p} (\prim_\ell H)^* \simeq \textnormal{gr}_F^{-p} \textrm{N}^\ell \prim_\ell(H^*) \simeq \textnormal{gr}_F^{-p+\ell} \prim_\ell(H^*).$$

\vspace{0.3cm}

\noindent
Consequently, duality translates as the transformation $(p,\ell) \mapsto (-p+\ell,\ell)$. Applying this rule, we deduce that \cite[Th.\,3(a)]{Fed18} is equivalent to

\vspace{-0.05cm}

$$\nu_{0,\mu_m,\ell}^p(H_{\bm{\alpha},\bm{\beta}})=\left\{
\begin{array}{cl}
  \ 1 \quad &\textnormal{if } p=-(\#\{i \ | \ \alpha_i \leq \alpha_m \}-\#\{j \ | \ \beta_j < \alpha_m \})+\ell_m(\bm{\alpha}) \\[1.5mm]
	&\textnormal{and } \ell=\ell_m(\bm{\alpha}) \\[1.5mm]
  \ 0 \quad &\textnormal{otherwise,}
\end{array}
\right.$$

\vspace{0.3cm}

\noindent
and this is equivalent to part (a) of the lemma.

\vspace{0.3cm}

\noindent
Similarly, \cite[Th.\,3(b)]{Fed18} is equivalent to

\vspace{-0.05cm}

$$\nu_{\infty,\lambda_m,\ell}^p(H_{\bm{\alpha},\bm{\beta}})=\left\{
\begin{array}{cl}
  \ 1 \quad &\textnormal{if } p=-(\#\{i \ | \ \alpha_i < \beta_m \}-\#\{j \ | \ \beta_j < \beta_m \})+\ell_m(\bm{\beta}) \\[2mm]
	&\textnormal{and } \ell=\ell_m(\bm{\beta}) \\[2mm]
  \ 0 \quad &\textnormal{otherwise,}
\end{array}
\right.$$

\vspace{0.3cm}

\noindent
and this is equivalent to part (b) of the lemma.
\end{proof}

\vspace{0.4cm}

It remains to show that the formulas of the previous lemma correspond to the formulas of Theorem \ref{fedorov}, up to a shift. This is a consequence of the following combinatorial lemma, insofar as $\#\{k \ | \ \alpha_{k} < \beta_{k} \}$ only depends on $\bm{\alpha}$ and $\bm{\beta}$.

\vspace{0.5cm}

\begin{lemm} We have the following relations:\\[1mm]
(i) $p(\bm{\alpha},\bm{\beta},\alpha_m)-(\#\{j \ | \ \beta_j < \alpha_m \} - \#\{i \ | \ \alpha_i < \alpha_m \})=\#\{k \ | \ \alpha_{k} < \beta_{k} \}$ \\[1mm]
(ii) $p(\bm{\alpha},\bm{\beta},\beta_m)-(\#\{j \ | \ \beta_j \leq \beta_m \} - \#\{i \ | \ \alpha_i < \beta_m \})=\#\{k \ | \ \alpha_{k} < \beta_{k} \}$.
\end{lemm}

\vspace{0.3cm}

\begin{proof} (i) Let us sum up in the following table the contributions of $k \in \{1,...,n\}$ to $p(\bm{\alpha},\bm{\beta},\alpha_m)$ and $\#\{j \ | \ \beta_j < \alpha_m \} - \#\{i \ | \ \alpha_i < \alpha_m \}$ according to the relative positions of $\alpha_k$, $\beta_k$ and $\alpha_m$.

\newpage

\begin{center}
\begin{tabular}{|c|c||c|c|} 
    \cline{3-4}
	  \multicolumn{2}{c|}{} & \multicolumn{2}{|c|}{contribution of $k$ to} \\
		\hline
    \multicolumn{2}{|c||}{relative positions} & \ $p(\bm{\alpha},\bm{\beta},\alpha_m)$ \ & \: $\#\{k \ | \ \beta_k < \alpha_m \}$ \: \\
		\multicolumn{2}{|c||}{} & & $- \#\{k \ | \ \alpha_k < \alpha_m \}$ \\
    \hline
		\hline
    \ $\alpha_k < \beta_k$ \ & \ \ $0 \leq \alpha_m < \alpha_k < \beta_k < 1$ \ \ & 1 & 0 \\
    \cline{2-4} 
        & \ \ $0 \leq \alpha_k = \alpha_m < \beta_k < 1$ \ \ & 1 & 0 \\
	  \cline{2-4}
				& \ \ $0 \leq \alpha_k < \alpha_m < \beta_k < 1$ \ \ & 0 & $-1$ \\
	  \cline{2-4}
				& \ \ $0 \leq \alpha_k < \beta_k < \alpha_m < 1$ \ \ & 1 & 0 \\
    \hline
		\hline
    \ $\alpha_k > \beta_k$ \ & \ \ $0 \leq \alpha_m < \beta_k < \alpha_k < 1$ \ \ & 0 & 0 \\
    \cline{2-4} 
        & \ \ $0 \leq \beta_k < \alpha_m < \alpha_k < 1$ \ \ & 1 & 1 \\
	  \cline{2-4}
				& \ \ $0 \leq \beta_k < \alpha_k = \alpha_m < 1$ \ \ & 1 & 1 \\
	  \cline{2-4}
				& \ \ $0 \leq \beta_k < \alpha_k < \alpha_m < 1$ \ \ & 0 & 0 \\
    \hline
\end{tabular}
\end{center}

\vspace{0.3cm}

\noindent
This table proves that $p(\bm{\alpha},\bm{\beta},\alpha_m)$ and $\#\{j \ | \ \beta_j < \alpha_m \} - \#\{i \ | \ \alpha_i < \alpha_m \}$ differ by $\#\{k \ | \ \alpha_{k} < \beta_{k} \}$, showing Formula (i).

\vspace{0.4cm}

\noindent
(ii) Let us now sum up in the following table the contributions of the integer $k$ to $p(\bm{\alpha},\bm{\beta},\beta_m)$ and $\#\{j \ | \ \beta_j \leq \beta_m \} - \#\{i \ | \ \alpha_i < \beta_m \}$ according to the relative positions of $\alpha_k$, $\beta_k$ and $\beta_m$.

\vspace{0.2cm}

\begin{center}
\begin{tabular}{|c|c||c|c|}
    \cline{3-4}
	  \multicolumn{2}{c|}{} & \multicolumn{2}{|c|}{contribution of $k$ to} \\
		\hline
    \multicolumn{2}{|c||}{relative positions} & \ $p(\bm{\alpha},\bm{\beta},\beta_m)$ \ & \: $\#\{k \ | \ \beta_k \leq \beta_m \}$ \: \\
		\multicolumn{2}{|c||}{} & & $- \#\{k \ | \ \alpha_k < \beta_m \}$ \\
    \hline
		\hline
    \ $\alpha_k < \beta_k$ \ & \ \ $0 \leq \beta_m < \alpha_k < \beta_k < 1$ \ \ & 1 & 0 \\
    \cline{2-4} 
        & \ \ $0 \leq \alpha_k < \beta_m < \beta_k < 1$ \ \ & 0 & $-1$ \\
	  \cline{2-4}
				& \ \ $0 \leq \alpha_k < \beta_k = \beta_m < 1$ \ \ & 1 & 0 \\
	  \cline{2-4}
				& \ \ $0 \leq \alpha_k < \beta_k < \beta_m < 1$ \ \ & 1 & 0 \\
    \hline
		\hline
    \ $\alpha_k > \beta_k$ \ & \ \ $0 \leq \beta_m < \beta_k < \alpha_k < 1$ \ \ & 0 & 0 \\
    \cline{2-4} 
        & \ \ $0 \leq \beta_k = \beta_m < \alpha_k < 1$ \ \ & 1 & 1 \\
	  \cline{2-4}
				& \ \ $0 \leq \beta_k < \beta_m < \alpha_k < 1$ \ \ & 1 & 1 \\
	  \cline{2-4}
				& \ \ $0 \leq \beta_k < \alpha_k < \beta_m < 1$ \ \ & 0 & 0 \\
    \hline
\end{tabular}
\end{center}

\vspace{0.2cm}

\noindent
This table proves that $p(\bm{\alpha},\bm{\beta},\beta_m)$ and $\#\{j \ | \ \beta_j \leq \beta_m \} - \#\{i \ | \ \alpha_i < \beta_m \}$ differ by $\#\{k \ | \ \alpha_{k} < \beta_{k} \}$, showing Formula (ii).
\end{proof}

\vspace{0.2cm}

\section*{Acknowledgements}

We thank first of all Claude Sabbah to whom this work owes a lot. The author is indebted as well to Michel Granger and Christian Sevenheck for their careful reading and constructive comments about this work. We also thank Michael Dettweiler for helpful discussions.

\bibliographystyle{amsalpha}
\bibliography{bibli}

\backmatter

\end{document}